
\input vanilla.sty
\scaletype{\magstep1}
\scalelinespacing{\magstep1}
\def\bull{\vrule height .9ex width .8ex depth -.1ex}


\pageno=-1

\title Common subspaces of $L_{p}$-spaces
\endtitle

\author Alexander Koldobsky\\Department of Mathematics \\
University of Missouri-Columbia\\Columbia, MO 65211
\endauthor

\vskip1truecm

\subheading{Abstract}  For $n\geq 2, p<2$ and $q>2,$ does there exist
an $n$-dimensional Banach space different from Hilbert spaces which
is isometric to subspaces of both $L_{p}$ and $L_{q}$? Generalizing the
construction from the paper "Zonoids whose polars are zonoids" by 
R.Schneider we give examples of such spaces. Moreover, for any compact
subset $Q$ of $(0,\infty)\setminus \{2k, k\in N\},$  we can construct a space  
isometric to subspaces of  $L_{q}$ for all $q\in Q$ simultaneously.

\vskip2truecm

\flushpar AMS classification: Primary 46B04, Secondary 46E30\newline
\newline e-mail: mathak\@mizzou1.missouri.edu 

\newpage
\pageno=1 

\title Common subspaces of $L_{p}$-spaces
\endtitle

\author Alexander Koldobsky\\Department of Mathematics \\
University of Missouri-Columbia\\Columbia, MO 65211
\endauthor

\vskip1truecm

\subheading{1. Introduction}

This work started with the following question:
For given $n\geq 2, p\in (0,2)$ and $q>2,$ does there exist an $n$-dimensional
Banach space which is different from Hilbert spaces and which is isometric
to subspaces of both $L_{p}$ and $L_{q}$?

It is a well-known fact first noticed by P.Levy that Hilbert spaces are 
isometric to subspaces of $L_{q}$ for all $q>0.$ On the other hand, it was 
proved in [4] that, for $n\geq 3, q>2, p>0,$ the function  
$\exp (-\|x\|_{q}^{p})$ is not positive definite where 
$ \|x\|_{q}=(|x_1|^q+...+|x_n|^q)^{1/q}.$ (This result gave an answer to 
a question posed by I.J.Schoenberg [10] in 1938) In 1966, J.Bretagnolle,
D.Dacunha-Castelle and J.L.Krivine [1] proved that,for $0<p<2$, 
a space $(E,\|\cdot \|)$ is isometric to a subspace of $L_p$ if and
only if the function $\exp (-\|x\|^{p})$ is positive definite. Thus, in
the language of isometries, the above mentioned result from [4] means
that, for every $n\geq 3, q>2, p\in (0,2),$ the space  $l_{n}^{q}$
is not isometric to a subspace of $L_p.$(For $p\geq 1,$ this fact was
first proved in [2]) The initial purpose of this work was to find a 
non-Hilbertian subspace $(E,\|\cdot \|)$ 
of $L_{q}$ with $q>2$ of the dimension at least 3 such that
the function $\exp (-\|x\|^{p})$ is positive definite. The latter problem 
is equivalent to that at the beginning of the paper. 

We prove, however, a more general fact: For every $n\geq 2$ and every compact 
subset $Q$ of $(0,\infty)\setminus \{2k, k\in N\},$  there exists 
an $n$-dimensional Banach space different from Hilbert spaces which is 
isometric to subspaces of  $L_{q}$ for all $q\in Q$ simultaneously.

In 1975, R.Schneider [9] proved that there exist  non-trivial zonoids 
whose polars are zonoids or, in other words, there exist non-Hilbertian
Banach spaces $X$ such that $X$ and $X^{*}$ are isometric to subspaces of $L_1.$
It turns out that Schneider's construction of special subspaces of $L_{1}$
can be extended to all numbers $q>0$ which are not even integers and 
in this way we obtain our main result.

\subheading{2. Some properties of spherical harmonics}

We start with some properties of spherical harmonics (see [6] for details).

Let $P_m$ denote the space of sperical
harmonics of degree $m$ on the unit sphere
$\Omega_n$ in $R^n.$ Remind that sperical
harmonics of degree $m$ are restrictions
to the sphere of harmonic homogeneous polynomials of degree $m.$ We consider
spherical harmonics as	functions from the space  $L_2(\Omega_n).$
Any two spherical harmonics of different degrees are orthogonal in 
$L_2(\Omega_n)$ [6, p.2]

The dimension $N(n,m)$ of the space $P_m$ can easily be calculated [6, p.4]:
$$N(n,m)={(2m+n-2)\Gamma(n+m-2)\over{\Gamma(m+1)\Gamma(n-1)}} \tag{1}$$

Let $\{Y_{mj}:j=1,...,N(n,m)\}$ be an orthonormal basis of the space $P_m.$
By the Addition Theorem [ 6, p.9], for every $x\in \Omega_n$,
$$ \sum_{j=1}^{N(n,m)} Y_{mj}^{2}(x)={N(n,m)\over \omega_n} \tag{2}$$
where $\omega_n={2\pi^{n/2}/\Gamma(n/2)}$ is the surface of the sphere
$\Omega_n.$

Linear combinations of functions $Y_{mj}$ are dense in the space
$L_2(\Omega_n)$ [ 6, p.43]. Therefore,
if $F$ is a continuous function on $\Omega_n$ and $(F,Y_{mj})=0$ for every $m=
0,1,2,...$ and every $j=1,....N(n,m)$ then $F\equiv 0$ on $\Omega_n.$ 
Here $( F,Y)$ stands for the scalar product in $L_2(\Omega_n).$

Let $\Delta$ be the Laplace-Beltrami operator on the sphere 
$\Omega_n.$ Then for every $Y_m\in P_m$ we have [ 6, p.39]   
$$\Delta Y_m + m(m+n-2)Y_m \equiv 0 \tag{3}$$

An immediate consequence of (3) ( and a well-known fact) is
 that $\Delta$ is a symmetric operator and we can
aplly Green's formula: for every function $H$ from the class $C^{2r},r\in N$
of functions on $\Omega_n$ having continuous partial derivatives of order $2r$
and for every $Y_m\in P_m,m\geq 1,$
$$ (-m(m+n-2))^{r} (H,Y_m)=(H,\Delta^{r} Y_m)=(\Delta^{r} H,Y_m) \tag{4}$$

We also need the Funk-Henke formula [ 6, p.20]:for every $Y_m\in P_m$,
every continuous function $f$ on $[-1,1]$ and every $x\in \Omega_n,$
$$\int_{\Omega_n}
f(\langle x,\xi \rangle) Y_{m} (\xi) d\xi=\lambda_{m} Y_{m}(x) \tag{5}$$
where $\langle \cdot,\cdot\rangle$ stands for the scalar product in $R^{n}$ and
$$\lambda_{m}={(-1)^{m} \pi^{(n-1)/2}\over 2^{m-1}\Gamma (m+(n-1)/2)}
\int_{-1}^{1} f(t) {d^m\over dt^m} (1-t^2)^{m+(n-3)/2} dt \tag{6}$$

Let us calculate $\lambda _m$ in the case where $f(t)=|t|^{q},q>0.$

\proclaim{Lemma 1} If $q>0,q\neq 2k,k\in N$ and  $f(t)=|t|^{q}$ then
$$ \lambda_m ={\pi^{n/2-1}\Gamma (q+1) \sin (\pi (m-q)/2) \Gamma ((m-q)/2)
\over 2^{q-1} \Gamma ((m+n+q)/2)} \tag{7} $$\endproclaim 

\demo{Proof} Assume first that $q>m$ and calculate the integral from (6)
by parts $m$ times. Then use the formula $\int_{-1}^{1} t^{2\alpha-1}(1-t^2)^
{\beta -1}dt
=\Gamma(\alpha)\Gamma(\beta)/\Gamma(\alpha+\beta)$ and formulas
for $\Gamma $-function: $\Gamma(2x)=2^{2x-1}\Gamma(x)\Gamma(x+1/2)/\pi^{1/2}$
and $\Gamma(1-x)\Gamma(x)=\pi/\sin(\pi x).$ We get (7) for $q>m.$ Note
that both sides of (7) are analytic functions of $q$ in the domain
$\Re q>0, q\neq 2k, k\in N.$ Because of the uniqueness of analytic extension,
(7) holds for every $q$ from this domain. We are done.\bull\enddemo

\subheading{3. Main result}

Let $X$ be an $n$-dimensional subspace of $L_q=L_{q}([0,1])$ with $q>0.$
Let $f_{1},...,f_{n}$ be a basis in $X$ and $\mu$ be the joint distribution
of the functions $f_{1},...,f_{n}$ with respect to Lebesgue measure ($\mu$
is a finite measure on $R^{n}).$ Then, for every $x\in R^{n}$,
$$\|x\|^{q}=\|\sum_{k=1}^{n} x_{k}f_{k}\|^{q}=\int_{0}^{1} |\sum_{k=1}^{n}
x_{k}f_{k}(t)|^{q} dt=\int_{R^{n}} |\langle x,\xi \rangle|^{q} d\mu (\xi)=
\int_{\Omega_{n}} |\langle x,\xi \rangle|^{q} d\nu (\xi) \tag{11}$$
where $\nu$ is the projection of $\mu$ to the sphere.(For every Borel subset
$A$ of $\Omega_{n}$, $\nu(A)=\int_{\{tA,t\in R\}} \|x\|_{2}^{q} d\mu(x) )$ 
The representation (11) of the norm is usially called the Levy representation.
It is clear  now that a norm in an $n$-dimensional Banach space admits the Levy
representation with a probability measure on the sphere if and only if 
this space is 
isometric to a subspace of $L_{q}.$ (Given the Levy representation we can
choose functions $f_{1},...,f_{n}$ on $[0,1]$ with the joint distribution 
$\nu$ and define an isometry by $x\rightarrow \sum_{k=1}^{n} x_{k}f_{k},
 x\in R^{n})$

If we replace the measure $\nu$ by an arbitrary continuous (not
necessarily non-negative) function on the sphere then a representation similar 
to the Levy representation is possible for a large  class of Banach spaces 
(see [5] for the Levy representation with distributions instead of measures;
such representation is possible for any Banach space and any $q$ which is not
an even integer) This is an idea going back to W.Blaschke that any smooth enough
function on the sphere can be represented in the form (11) with a continuous
function instead of a measure on the sphere. However, W.Blaschke and then
R.Schneider [8] restricted themselves to the case $q=1$ which is particularly
important in the theory of convex bodies. The following theorem is an
extension of R.Schneider's results from [ 8, p.77] and [ 9, p.367] to all 
positive numbers $q$ which are not even integers 

\proclaim{Theorem 1} Let $q>0, q\neq 2k, k\in N$
and let $H$ be an even function of
the class $C^{2r}$ on $\Omega_{n}$ where $r\in N$ and  $2r>n+q.$
Then there exists a continuous function $b_{H}$ on the
sphere  $\Omega_{n}$ such that, for every $x\in \Omega_{n},$
$$H(x)=\int_{\Omega_{n}} |\langle x,\xi \rangle|^{q}b_{H}(\xi)d\xi \tag{8}$$
Besides that, there exist constants $K(q)$ and $L(q)$ depending on $n$ and $q$
only such that, for every  $x\in \Omega_{n},$
$$|b_{H}(x)|\leq K(q)\|H\|_{L_{2}(\Omega_{n})} + L(q)\|\Delta^{2r}H\|_{L_{2}
(\Omega_{n})} \tag{9}$$
\endproclaim

\demo{Proof} Define a function $b_{H}$ on $\Omega_{n}$ by
$$b_{H}(x)=\sum_{m=0}^{\infty} \lambda_{m}^{-1}\sum_{j=1}^{N(n,m)} (N,Y_{mj})
Y_{mj}(x) \tag{10}$$
Since $H$ is an even function and $Y_{mj}$ are odd functions if $m$ is odd
the sum is, actually, taken over even integers $m$ only.

Let us prove that the series in the
right-hand side of (10) converges uniformly on
$\Omega_{n}.$ By the Cauchy-Schwartz inequality, (2) and
the fact that $Y_{mj}$ form an orthonormal basis in $P_{m}$, we get
$$|\sum_{j=1}^{N(n,m)} (\Delta^{r}N,Y_{mj})Y_{mj}(x)|\leq 
(\sum_{j=1}^{N(n,m)} (\Delta^{r}N,Y_{mj})^2)^{1/2} (\sum_{j=1}^{N(n,m)}
Y_{mj}^{2}(x))^{1/2}$$
$$\leq \|\Delta^{r}H\|_{L_{2}(\Omega_{n})}\bigl({N(n,m)\over
\omega_n}\bigr)^{1/2}$$
It follows from (4) and the latter inequality that
$$|b_{H}(x)|\leq |\lambda_{0}^{-1}(H,Y_{0})Y_{0}(x)|+\sum_{m=2;2|m}^{\infty}
\lambda_{m}^{-1}({-1\over m(m+n-2)})^{r}|\sum_{j=1}^{N(n,m)}
(\Delta^{r}H,Y_{mj})
Y_{mj}(x)| \leq$$
$$ |\lambda_{0}|^{-1} \omega_{n}^{-1/2}\|H\|_{L_{2}(\Omega_{n})}+
\sum_{m=2;2|m}^{\infty}\lambda_{m}^{-1}
m^{-2r}\bigl({N(n,m)\over\omega_n}\bigr)^{1/2}
\|\Delta^{r}H\|_{L_{2}(\Omega_{n})}$$  
Let us show that the series 
$\sum_{m=2;2|m}^{\infty}\lambda_{m}^{-1} m^{-2r}({N(n,m)/\omega_n})^{1/2}$
converges. In fact, it follows from (1) that $N(n,m)=O(m^{n-2})$ and 
it is an easy consequence of (7) and the Stirling formula that
 $\lambda_{m}^{-1}=O(m^{(n+2q)/2}).$
Since $2r>n+q=(n+2q)/2 + (n-1)/2 + 1$ we get   	 
$\lambda_{m}^{-1} m^{-2r}({N(n,m)/\omega_n})^{1/2}=o(m^{-1-\epsilon})$ for
some $\epsilon >0$, and the series is convergent. We denote the sum
of this series by $L(q)$ and put $K(q)=|\lambda_{0}|^{-1} \omega_{n}^{-1/2}$, 
so we get (9).

We have proved that the series in (10) converges uniformly and defines
a continuous function 
on $\Omega_{n}.$ It follows from (5) and the fact that all functions $Y_{mj}$ 
are orthogonal that $(H,Y_{mj})=(\int_{\Omega_{n}} |\langle x,\xi \rangle|^{q}
b_{H}(\xi)d\xi, Y_{mj}(x))$ for every $m=0,1,2,...$ and $j=1,...,N(n,m).$ Hence,
the function $b_{H}$ satisfies (8). \bull \enddemo

Let $X$ be an $n$-dimensional Banach space, $q>0, q\neq 2k, k\in N.$ 
Let $c(q)=\Gamma ((n+q)/2)/(2\Gamma((q+1)/2) 
\pi^{(n-1)/2})$ be a constant such that
$1=c(q)\int_{\Omega_{n}}|\langle x,\xi \rangle|^{q} d\xi$ for every
$x\in \Omega_{n}.$ (The latter integral
does not depend on the choice of $x\in \Omega_{n}$; it means that the norm
of the space $l_{2}^{n}$ admits the Levy representation with the uniform 
measure on the sphere and the space $l_{2}^{n}$ is isometric to a subspace 
of $L_{q}$ for every $q$)

Denote by $H(x), x\in\Omega_{n}$ the restriction of the function
 $\|x\|^{q}$ to the sphere $\Omega_{n}$.
Assume that the function $H$ belongs to the class $C^{2r}$ on $\Omega_{n}$
where $2r>n+q, r\in N.$  
Let $b_{H}$ be the function corresponding to $H$ by Theorem 1.

\proclaim{Lemma 2} If the number $K(q)\|H-1\|_{L_{2}(\Omega_{n})} +
L(q)\|\Delta^{r}H\|_{L_{2}(\Omega_{n})}$ is less than $c(q)$ then
the space $X$ is isometric to a subspace of $L_{q}.$ \endproclaim

\demo{Proof} By (8) and definition of the number $c(q)$,
$$H(x)-1=\int_{\Omega_{n}} |\langle x,\xi \rangle|^{q}(b_{H}(\xi)-c(q))d\xi$$
for every $x\in \Omega_{n}.$ By (9), $|b_{H}(x)-c(q)| < c(q)$ for every 
$x\in \Omega_{n}.$ It means that the function $b_{H}$ is positive on the
sphere. The equality (8) means that the space $X$ admits the 
Levy representation with a non-negative measure and, by the reasoning at
the beginning of Section 3, $X$ is isometric to a subspace of $L_{q}.$
\bull \enddemo

Now we are able to prove the main result of this paper. Let us only note
that, for every function $f$ of the class $C^{2}$ on the sphere  
$\Omega_{n}$ and for a small enough number $\lambda,$ the function
$N(x)=1 + \lambda f(x), x\in \Omega_{n}$ is the restriction to the sphere
of some norm in $R^{n}.$ This is an easy consequence of the following
one-dimensional fact: If $a,b\in R,$ $g$ is a convex function on $[a,b]$
with $g''>\delta >0$ on $[a,b]$ for some $\delta$ and $h\in C^{2}[a,b]$ then
functions $g + \lambda h$ have positive second derivatives on $[a,b]$
for sufficiently small $\lambda$'s, and, hence, are convex on $[a,b].$ 

\proclaim{Theorem 2} Let $Q$ be a compact subset of $(0,\infty)\setminus
\{2k, k\in N\}.$ Then there exists a Banach space different from
Hilbert spaces which is isometric to a subspace of $L_{q}$ for every $q\in Q.$
\endproclaim

\demo{Proof} Let $f$ be any infinitely differentiable function on $\Omega_{n}$
and fix a number $r\in N$ so that $2r > n+q$  for every $q\in Q.$
Choose a sufficiently small number $\lambda$ such that the function
$N(x)=1 + \lambda f(x), x\in \Omega_{n}$ is the restriction to the sphere
of some norm in $R^{n}$ ( see the remark before Theorem 2) and such that,
for every $q\in Q,$
the function $H(x)=(N(x))^{q}$ satisfies the condition of Lemma 2.     
The possibility of such choice of $\lambda$ follows from the facts that
$K(q), L(q)$ and $c(q)$ are continuous functions of $q$ on the set $Q$ and
that $\|H-1\|_{L_{2}(\Omega_{n})}$ and $\|\Delta^{2r}H\|_{L_{2}(\Omega_{n})}$
tend to zero uniformly with resrect to $q\in Q$ as $\lambda$ tends to zero.
Now we can apply Lemma 2 to complete the proof. \bull \enddemo

Finally, let us consider the case where $q$ is an even integer. It is
easy to see that, for any fixed number $2k, k\in N, k>1,$ we can make the space
$X$ constructed in Theorem 2 isometric to a subspace of $L_{2k}.$ In fact,
let $N(x)=(1 + \lambda (x_{1}^{2k}+...+x_{n}^{2k}))^{1/4}.$ For sufficiently 
small numbers $\lambda,$ $N$ is the restriction to the sphere of some norm 
in $R^{n}$ and the corresponding space $X$ is isometric to a subspace
of $L_{q}$ for every $q\in Q.$ On the other hand, $X$ is isometric to a subspace
of $L_{2k}$ because the norm admits the Levy representation with a measure on
the sphere:
$$ 1 + \lambda (x_{1}^{2k}+...+x_{n}^{2k}) = \int_{\Omega_{n}}
 |\langle x,\xi \rangle|^{2k}(c(2k)d\xi +\lambda d\delta_{1}(\xi) +...+
\lambda d\delta_{n}(\xi))$$
where $\delta_{i}$ is a unit mass at the point $\xi \in R_{n}$ with
$\xi_{i}=1, \xi_{j}=0, j\neq i.$

Let us show that one can not make the space $X$ isometric to subspaces 
of $L_{2p}$ and $L_{2q}$ if $p,q\in N$ and do not
have common factors. In fact, if $(X, \|\cdot\|)$ is such a space then,
for every $x\in R_{n},$
$$ \|x\|^{4pq} = (\int_{\Omega_{n}} |\langle x,\xi \rangle|^{2p} d\mu(\xi))^{2q}
=(\int_{\Omega_{n}} |\langle x,\xi \rangle|^{2q} d\nu(\xi))^{2p}$$
for some measures $\mu, \nu$ on $\Omega_{n}.$ The functions in the latter 
equality are polynomials and, since the polynomial ring has the unique 
factorization property , we conclude that $\|x\|^{2}$ is a homogeneous 
polynomial of the second order and $X$ is a Hilbert space.

The situation is not clear if $p$ and $q$ have common factors. One can find 
some interesting results on Banach spaces with polynomial norms and 
on the structure of subspaces of $L_{2k}, k\in N$ in the paper [7].   

\subheading{Acknowledgements}
I wish to thank Prof. Nigel Kalton for valuable remarks and
helpful discussions during the work on this problem. I am grateful to
Prof. Hermann Konig for bringing the paper [9] to my attention. 

\vskip2truecm

\newpage
\subheading{References} 

\item{1. } Bretagnolle, J., Dacunha-Castelle, D. and Krivine, J.L.:
{\it Lois stables et espaces $L_{p},$} Ann. Inst. H.Poincare, Ser.B
   2(1966), 231-259.

\item{2. } Dor, L.: {\it Potentials and isometric embeddings in $L_{1},$}
Israel J. Math.  24 (1976), 260-268.

\item{3. } Grzaslewicz, R.: {\it Plane sections of the unit ball of $L_{p},$}
Acta Math. Hung. 52(1988), 219-225.

\item{4. } Koldobsky, A.: {\it Schoenberg's problem on positive definite 
functions,} Algebra and Analysis 3 (1991), No. 3, 78-85 (Russian);
English translation in  St.Petersburg Math. J. 3 (1992), 563-570

\item{5. } Koldobsky, A.: {\it Generalized Levy representation of norms and
isometric	embeddings into Lp-spaces,} Ann. Inst. H.Poincare (Prob. and Stat.)
28 (1992), 335-353

\item{6. } Muller, C.: {\it Spherical Harmonics,} Lect. Notes in Math.
17 , Springer-Verlag, Berlin, 1966.

\item{7. } Reznick, B.: {\it Banach spaces with polynomial norms,}
Pacific J. Math. 82 (1979), 223-235.

\item{8. } Schneider, R.: {\it Zu einem problem von Shephard uber die
projektionen konvexer korper,} Math. Zeitschrift 101 (1967), 71-82.

\item{9. } Schneider, R.:  {\it Zonoids whose polars are zonoids,}
Proc. Amer. Math. Soc. 50 (1975), 365-368.

\item{10. } Schoenberg, I.J.: {\it Metric spaces and positive definite
functions,} Trans. Amer. Math. Soc. 44 (1938), 522-536.

\bye